\documentclass[11pt]{amsart}
\usepackage{amssymb}
\usepackage{epsfig}
\usepackage{graphicx}

\usepackage{fullpage}
\usepackage[numbers,sort&compress]{natbib}
\usepackage{color}
\usepackage{setspace}
\usepackage[marginal]{footmisc}

\begin{document}
\newtheorem*{pro*}{Proposition~$3.7^{(')}$}
\newtheorem{mthm}{Theorem}
\newtheorem{mcor}{Corollary}
\newtheorem{mpro}{Proposition}
\newtheorem{mfig}{figure}
\newtheorem{mlem}{Lemma}
\newtheorem{mdef}{Definition}
\newtheorem{mrem}{Remark}
\newtheorem{mpic}{Picture}
\newtheorem{rem}{Remark}[section]
\newcommand{\ra}{{\mbox{$\rightarrow$}}}
\newtheorem{Remark}{Remark}[section]
\newtheorem{thm}{Theorem}[section]
\newtheorem{pro}{Proposition}[section]
\newtheorem*{proA}{Proposition A}

\newtheorem{lem}{Lemma}[section]
\newtheorem{defi}{Definition}[section]
\newtheorem{cor}{Corollary}[section]

\title[Nonexistence of anti-symmetric solutions]{Nonexistence of anti-symmetric solutions for fractional Hardy-H\'{e}non System}
\author[J. Hu and Z. Du]{Jiaqi Hu$^\dag$ and Zhuoran Du$^\ddag$}

\thanks{$^\dag$ School of Mathematics, Hunan University, Changsha 410082, China.
{Email: hujq@hnu.edu.cn}.}
\thanks{}

\thanks{$^\ddag$School of Mathematics, Hunan University, Changsha 410082, China.
{Email: duzr@hnu.edu.cn}. }
\thanks{}

\date{\today}

\thanks{}

\date{}

\maketitle

\begin{abstract}
We study  anti-symmetric solutions about the hyperplane $\{x_n=0\}$ to the following fractional Hardy-H\'{e}non system
$$
\left\{\begin{aligned}
&(-\Delta)^{s_1}u(x)=|x|^\alpha v^p(x),\ \ x\in\mathbb{R}_+^n,
\\&(-\Delta)^{s_2}v(x)=|x|^\beta u^q(x),\ \ x\in\mathbb{R}_+^n,
\\&u(x)\geq 0,\ \ v(x)\geq 0,\ \ x\in\mathbb{R}_+^n,
\end{aligned}\right.
$$
where $0<s_1,s_2<1$, $n>2\max\{s_1,s_2\}$. Nonexistence of anti-symmetric solutions are obtained in some appropriate domains of $(p,q)$ under some corresponding assumptions of $\alpha,\beta$ via the methods of moving spheres and moving planes. Particularly, for the case $s_1=s_2$, one of our results shows that one domain of $(p,q)$, where nonexistence of anti-symmetric solutions with  appropriate decay conditions holds true, locates at above the fractional Sobolev's hyperbola under appropriate condition of $\alpha, \beta$.

\end{abstract}


\section{Introduction}

In this paper, we study anti-symmetric solutions  about the hyperplane $\{x_n=0\}$ to the following system involving fractional Laplacian
\begin{equation}\label{SI3}
\left\{\begin{aligned}
&(-\Delta)^{s_1}u(x)=|x|^\alpha v^p(x),\ \ x\in\mathbb{R}_+^n,
\\&(-\Delta)^{s_2}v(x)=|x|^\beta u^q(x),\ \ x\in\mathbb{R}_+^n,
\\&u(x)\geq 0,\ \ v(x)\geq 0,\ \ x\in\mathbb{R}_+^n,
\\&u(x',x_n)=-u(x',-x_n),\ \ v(x',x_n)=-v(x',-x_n),\ \ \ \ x=(x',x_n)\in\mathbb{R}^n,
\end{aligned}\right.
\end{equation}
where $s_1,s_2\in(0,1)$, $n>\max\{2s_1,2s_2\}$, $\mathbb{R}_+^n=\{x=(x',x_n)\in\mathbb{R}^n|x_n>0\}$ and $x'=(x_1,x_2,\cdots,x_{n-1})$.

The fractional Laplacian $(-\Delta)^s \, (0<s<1)$  is a nonlocal operator defined by
$$(-\Delta)^su(x)=C(n,s)P.V.\int_{\mathbb{R}^n}\frac{u(x)-u(y)}{|x-y|^{n+2s}}dy,$$
where $P.V.$ stands for the Cauchy principal value and $C(n,s)=\Big(\int_{\mathbb{R}^n}\frac{1-\cos\xi}{|\xi|^{n+2s}}d\xi\Big)^{-1}$ (see \cite{Caff, H}).
Let
$$L_{2s}=\Big\{u:\mathbb{R}^n\to\mathbb{R}|\int_{\mathbb{R}^n}\frac{|u(x)|}{1+|x|^{n+2s}}dx<+\infty\Big\}.$$
Then for $u\in L_{2s}, (-\Delta)^su$ can be defined in distributional sense (see \cite{S})
$$\int_{\mathbb{R}^n}(-\Delta)^su \varphi dx=\int_{\mathbb{R}^n}u(-\Delta)^s\varphi dx,\ \ \ \ \text{for any}\ \ \varphi\in \mathcal{S}.$$
Moreover, $(-\Delta)^su$ is well defined for $u\in L_{2s}\cap C_{loc}^{1,1}(\mathbb{R}^n)$.
We call $(u,v)$  a classical solution of (\ref{SI3}) if $(u,v)\in( L_{2s_1}\cap  C^{1,1}_{loc}(\mathbb{R}_+^n)\cap C(\mathbb{R}^n)) \times(L_{2s_2}\cap  C^{1,1}_{loc}(\mathbb{R}_+^n)\cap C(\mathbb{R}^n))$ and satisfies (\ref{SI3}).

As well known, the method of moving planes and moving spheres play an important role in proving the nonexistence of solutions. Chen et al. \cite{Chen3, Chen4} introduced direct method of moving planes and moving spheres for fractional Laplacian, which have been widely applied to derive the symmetry, monotonicity and nonexistence and even a prior estimates of solutions for some equations involving fractional Laplacian. In such process, some suitable forms of maximum principles are the key ingredients. The method of moving planes in integral forms is also a vital tool for classification of solutions (see \cite{Chen5}).

Recently, Li and Zhuo in \cite{LZ1} classified anti-symmetric classical solutions of Lane-Emden system (\ref{SI3}) in the case of $s_1=s_2=:s\in (0,1)$ and $\alpha=\beta=0$. They established the following Liouville type theorem.
\begin{pro}\label{P01}\cite{LZ1}
Given $0<p,q\leq\frac{n+2s}{n-2s},$  assume that $(u, v)$ is an anti-symmetric classical solution of
system (\ref{SI3}). If $0 < pq < 1$ or $p + 2s > 1$ and $q + 2s > 1$, then $(u,v)\equiv (0,0)$.
\end{pro}
\noindent As a corollary of Proposition 1, the nonexistence results in the larger space $L_{2s+1}$ follows immediately for the case $p+2s>1$, $q+2s>1$.
\begin{pro}\label{P02}\cite{LZ1}
Assume that $u$ and $v \in L_{2s+1}\cap C^{1,1}_{loc}(\mathbb{R}^n_+)\cap C(\mathbb{R}^n)$ satisfy the system (\ref{SI3}). Then if $0 < p, q \leq\frac{n+2s}{n-2s}$,  $ p+2s > 1$ and $q + 2s > 1$, $(u,v)\equiv (0,0)$ is the only solution.
\end{pro}
The nonexistence  of anti-symmetric classical solutions to the corresponding scalar problem  were shown in \cite{Zhuo}.

The following Hardy-H\'{e}non system with homogeneous Dirichlet boundary conditions has been investigated widely
\begin{equation}\label{I1}
\left\{\begin{aligned}
&(-\Delta)^{s_1}u(x)=|x|^\alpha u^p, \ \ u(x)\geq 0,\ \ &x\in\Omega,
\\&(-\Delta)^{s_2}v(x)=|x|^\beta v^q, \ \ v(x)\geq 0,\ \ &x\in\Omega,
\\&u(x)=v(x)=0,\ \ \ \ &x\in\mathbb{R}^n\setminus\Omega.
\end{aligned}\right.
\end{equation}
There are enormous nonexistence results of  (\ref{I1}) for the case $\Omega=\mathbb{R}^n$. We list some main results as follows.

If $s_1=s_2=1$, for $\alpha,\beta\geq0$, system (\ref{I1}) is the well-known H\'{e}non-Lane-Emden system. It has been conjectured that the Sobolev's hyperbola
$$\Big\{p>0,q>0:\frac{n+\alpha}{p+1}+\frac{n+\beta}{q+1}=n-2\Big\}$$
is the critical dividing curve between existence and nonexistence of solutions to (\ref{I1}). Particularly, the H\'{e}non-Lane-Emden conjecture states that the system (\ref{I1}) admits no nonnegative non-trivial solutions if $p>0$, $q>0$ and $\frac{n+\alpha}{p+1}+\frac{n+\beta}{q+1}>n-2$. For $\alpha=\beta=0$, this conjecture has been completely proved for radial solutions (see \cite{M,SZ1}). However, for non-radial solutions, the conjecture is only fully answered when $n\leq 4$ (see \cite{ PO, SZ, SP}). In higher dimensions, the conjecture was partially solved. Felmer and Figueriredo \cite{FF} showed that the system (\ref{I1}) admits no classical positive solutions if
$$0<p,q\leq \frac{n+2}{n-2}\ \ \ \ \text{and  }\ \ (p,q)\neq\Big(\frac{n+2}{n-2},\frac{n+2}{n-2}\Big).$$
Busca and Man\'{a}sevich \cite{BM} proved  the conjecture  if
$$\alpha_1,\alpha_2\geq \frac{n-2}{2}\ \  \text {and}\ \ (\alpha_1,\alpha_2)\neq \Big(\frac{n-2}{2},\frac{n-2}{2}\Big),$$
where
$$\alpha_1=\frac{2(p+1)}{pq-1},\ \ \ \ \alpha_2=\frac{2(q+1)}{pq-1}, \ \  pq>1.$$
When $\alpha,\beta>0$, Fazly and Ghoussoub \cite{FG} showed that the conjecture hold for dimension $n = 3$ under the assumption of the
boundness of positive solutions, Li and Zhang \cite{KL} removed this assumption and
proved this conjecture for dimension $n = 3$. When $\min\{\alpha,\beta\}>-2$, the conjecture is proved for bounded solutions in $n=3$ (see \cite{P}).

If $s_1=s_2=:s\in(0,1), \alpha, \beta\geq0$, there are fewer nonexistence results of solutions to  system (\ref{I1}) in the case of $p>0$, $q>0$ and $\frac{n+\alpha}{p+1}+\frac{n+\beta}{q+1}>n-2s$, namely the case that $(p,q)$ locates at bottom left of the fractional Sobolev's hyperbola
$\Big\{p>0,q>0:\frac{n+\alpha}{p+1}+\frac{n+\beta}{q+1}=n-2s\Big\}$. For $\alpha=\beta=0$,  Quass and Xia in \cite{QX1} proved that there exist no classical positive solutions to  (\ref{I1}) provided that
\begin{equation}\label{SE1}
\alpha^s_1,\alpha^s_2\in\Big[\frac{n-2s_1}{2},n-2s_1\Big), \ \ \ \ \text{and }\ \ (\alpha^s_1,\alpha^s_2)\neq \Big(\frac{n-2s}{2},\frac{n-2s}{2}\Big),
\end{equation}
where $$\alpha_1^s=\frac{2s(q+1)}{pq-1}, \ \ \alpha_2^s=\frac{2s(p+1)}{pq-1},\ \ p,q>0,\ \  pq>1.$$
Note that the region (\ref{SE1}) of $(p,q)$ contains the following region
$$\Big\{(p,q):\frac{n}{n-2s}<p,q\leq \frac{n+2s}{n-2s},\ \ \text{and }\ \ (p,q)\neq \Big(\frac{n+2s}{n-2s},\frac{n+2s}{n-2s}\Big)\Big\}.$$
As $\min\{\alpha,\beta\}>-2s$, Peng \cite{PS1} derived that the system ({\ref{I1}) admit no nonnegative classical solutions if $0<p<\frac{n+2s+2\alpha}{n-2s}$ and $0<q<\frac{n+2s+2\beta}{n-2s}$.

For scalar equation (i.e., $s_1=s_2:=s, \alpha=\beta, p=q, u=v$), in the Laplacian case,  if $\alpha=0$, a celebrated Liouville type theorem was showed by Gidas and Spruck \cite{GS} for $1<p<\frac{n+2}{n-2}$;  if $\alpha\leq -2$ and $p>1$, there is no any positive solution (see \cite{G,N}); if $\alpha>-2$ and $1<p<\frac{n+2+2\alpha}{n-2}$, Phan and Souplet \cite{PS} derived a Liouville theorem for bounded solutions; if  $0 <p \leq 1$, the nonexistence result was proved by Dai and Qin \cite{DQ1} for any $\alpha$. We also refer to \cite{GB,ME} and references therein. For the scalar fractional Laplacian case, Chen, Li and Li \cite{Chen3}, Jin, Xiong and Li \cite{J} proved the nonexistence results for $\alpha=0$ and $0<p<\frac{n+ 2s}{n-2s}.$ If $\alpha>-2s$, Dai and Qin \cite{DQ}  showed a Liouville type theorem for optimal range $0<p<\frac{n+2s+2\alpha}{n-2s}$.

For the case $\Omega=\mathbb{R}^n_+$, there are the following several main results.

If $s_1=s_2=1$, $\alpha=\beta=0$, $\min\{p,q\}>1$,  a Liouville Theorem is proved for bounded solutions by Chen, Lin and Zou \cite{CL}. If $s_1=s_2=:s \in(0,1)$, for $\alpha=\beta=0$, the nonexistence of  positive viscosity bounded solutions to system (\ref{I1}) was showed by Quaas and Xia in \cite{QX}. For $\alpha,\beta>-2s$, if $p\geq \frac{n+2s+\alpha}{n-2s}$ and $q\geq \frac{n+2s+\beta}{n-2s}$, Duong and Le \cite{DL} obtained the nonexistence of solutions satisfying the following decay at infinity
$$u(x) = \mbox{o}(|x|^{-\frac{4s+\beta}{q-1}}) \ \ \text{and} \ \ v(x) = \mbox{o}(|x|^{-\frac{4s+\alpha}{p-1}}).$$  For general $s_1, s_2\in(0,1)$, Le in \cite{L} concluded a Liouville type theorem. Precisely, they obtained that if $1\leq p\leq \frac{n+2s_1+2\alpha}{n-2s_2}$, $1\leq q\leq\frac{n+2s_2+2\alpha}{n-2s_1}$ and $(p,q)\neq (\frac{n+2s_1+2\alpha}{n-2s_2},\frac{n+2s_2+2\alpha}{n-2s_1}$), $\alpha>-2s_1$ and $\beta>-2s_2$, then $(u,v)\equiv (0,0)$ is the only nonnegative classical solutions to system (\ref{I1}). More nonexistence results  for general nonlinearities in a half space can be seen in \cite{DP,ZL}.
For the corresponding scalar problem of (\ref{I1}) with Laplacian and $\alpha=\beta=0$, Gidas and Spruck  \cite{GS} obtained the nonexistence of nontrivial nonnegative classical solution of (\ref{I1}) for $1<p\leq\frac{n+2}{n-2}$. For  the corresponding scalar problem  (\ref{I1}) with fractional Laplacian and $\alpha=\beta=0$, Chen, Li and Li \cite{Chen3} showed that $u\equiv 0$ is the only nonnegative solution to (\ref{I1}) for $1<p\leq\frac{n+2s}{n-2s}$.  Recently a Liouville type theorem of the corresponding scalar problem for $1<p<\frac{n+2s+2\alpha}{n-2s}, \alpha>-2s$ and $s\in (0,1]$ was established by Dai and Qin in \cite{DQ}.

In this paper we will study nonexistence of anti-symmetric classical solutions to  the system (\ref{SI3}) for general $\alpha, \beta, p, q$.

For $\alpha>-2s_1$ and $\beta>-2s_2$, we denote
$$
\mathcal{R}_{sub}:=\Big\{(p,q)|0<p\leq \frac{n+2s_1+2\alpha}{n-2s_2}, 0<q\leq \frac{n+2s_2+2\beta}{n-2s_1}, (p,q)\neq \Big(\frac{n+2s_1+2\alpha}{n-2s_2},\frac{n+2s_2+2\beta}{n-2s_1}\Big)\Big\}.
$$
Note that for the case $s_1=s_2$, the set $\mathcal{R}_{sub}$ locates at bottom left of the preceding fractional Sobolev's hyperbola.

Throughout this paper, we always assume $s_1,s_2\in(0,1)$, $n>\max\{2s_1,2s_2\}$.
and use $C$ to denote a general positive constant whose value may vary from line to line even the same line.
 Our main results are as follows.

\begin{thm}\label{SST1}
For $(p,q)\in \mathcal{R}_{sub}$, assume that $(u,v)$ is a  classical solution of the system (\ref{SI3}). For either one of the following two cases\\
(\romannumeral 1) $\min\{p+2s_1,p+2s_1+\alpha\}>1$ and $\min\{q+2s_2,q+2s_2+\beta\}>1$, \\
(\romannumeral 2) $0<pq<1$ with $\alpha\geq-2s_1pq$, $\beta\geq-2s_2pq$,\\
we have that $(u,v)=(0,0)$.
\end{thm}
The nonexistence results (i) of Theorem 1 can be extended to a larger space.
\begin{thm}\label{T01}
Assume that  $(p,q)\in \mathcal{R}_{sub}$ and $(u,v)\in (L_{2s_1+1}\cap C^{1,1}_{loc}(\mathbb{R}^n_+)\cap C(\mathbb{R}^n))\times (L_{2s_2+1}\cap C^{1,1}_{loc}(\mathbb{R}^n_+)\cap C(\mathbb{R}^n))$ satisfies the system (\ref{SI3}). Then for the case that $\min\{p+2s_1,p+2s_1+\alpha\}>1$ and $\min\{q+2s_2,q+2s_2+\beta\}>1$, $(u,v)\equiv (0,0)$ is the only solution.
\end{thm}
Combining our anti-symmetric property, in the process of Theorem \ref{SST1}, we only utilized the extended spaces $L_{2s_1+1}, L_{2s_2+1}$ instead of the usual spaces $L_{2s_1}, L_{2s_2}$ in the case that $\min\{p+2s_1,p+2s_1+\alpha\}>1$ and $\min\{q+2s_2,q+2s_2+\beta\}>1$. One can see that Theorem \ref{T01} is a direct corollary of (i) of Theorem \ref{SST1}.
\begin{rem}
Our results of Theorem \ref{SST1} and Theorem \ref{T01} are the extension to general $s_1,s_2,\alpha,\beta$ of the nonexistence results of Li and Zhuo \cite{LZ1} (see preceding Propositions \ref{P01} and  \ref{P02}) except one critical point of $(p,q)$.
\end{rem}

\begin{rem}
 When $s_1=s_2$, $p=q$, $\alpha=\beta$ and $u=v$, the results of Theorems \ref{SST1} and  \ref{T01} are the nonexistence of nontrivial  classical solutions to the corresponding scalar problem.
\end{rem}

Under appropriate decay conditions of $u$ and $v$ at infinity, we can extend the nonexistence result of classical solutions of (\ref{SI3}) to an unbounded domain of $(p,q)$. Particularly, this unbounded domain, except at most a bounded sub-domain,  locates at above the preceding fractional Sobolev's hyperbola for the case $s_1=s_2$.
\begin{thm}\label{T2}
Suppose $p\geq \frac{n+2s_1+\alpha}{n-2s_2}$, $q\geq \frac{n+2s_2+\beta}{n-2s_1}$, $\alpha>-2s_2$, $\beta>-2s_1$. Assume $(u,v)$ is a  classical solution of system (\ref{SI3}) satisfying
$$\mathop{\overline{\lim}}\limits_{x\to \infty}\frac{u(x)}{|x|^{a}}\leq C\ \  \text{and}\ \ \mathop{\overline{\lim}}\limits_{x\to \infty}\frac{v(x)}{|x|^{b}}\leq C,$$
for some $C>0$, where $a=-\frac{2s_1+2s_2+\beta}{q-1}$ and $b=-\frac{2s_1+2s_2+\alpha}{p-1}$. Then $(u,v)\equiv(0,0)$.
\end{thm}

\begin{rem}
The results of Theorem \ref{T2} are new even if for the corresponding scalar problem with $\alpha=0$.
\end{rem}
\begin{rem}
When $\alpha,\beta$ are positive, define the region
\begin{eqnarray}\nonumber
\mathcal{R}_{sup}:=&&\left\{(p,q)|\frac{n+2s_1+\alpha}{n-2s_2}\leq p\leq \frac{n+2s_1+2\alpha}{n-2s_2},\frac{n+2s_2+\beta}{n-2s_1}\leq p\leq\frac{n+2s_2+2\beta}{n-2s_1}, \right.\\ \nonumber
&&\left.
 (p,q)\neq \left(\frac{n+2s_1+2\alpha}{n-2s_2},\frac{n+2s_2+2\beta}{n-2s_1}\right)\right\}.
\end{eqnarray}
Note that $\mathcal{R}_{sup}$ is contained in the nonexistence region of $(p,q)$ obtained in Theorem \ref{SST1}. Hence, if $(p,q)\in \mathcal{R}_{sup}$, Theorem \ref{SST1} tells us that the results of Theorem \ref{T2} still hold true without the decay conditions.
\end{rem}

\section{Preliminaries}
In this section, we introduce and prove some necessary lemmas.

\begin{pro}(\cite{DL})\label{L5}
Let  $s\in(0,1)$ and $w(y)\in L_{2s}\cap C_{loc}^{1,1}(\mathbb{R}^n)$ satisfy
$w(y)=-w(y^-_\lambda)$, where $y^-_\lambda=(y',
2\lambda-y_n)$ for any real number $\lambda$. Assume there exists $x\in \Sigma_\lambda$ such that
$$w(x)=\mathop{\inf}\limits_{\Sigma_\lambda }w(y)<0 \ \ \text{ and } \ \
(-\Delta)^sw(x)+c(x)w(x)\geq 0,$$
where $\Sigma_\lambda=\{x\in\mathbb{R}^n|x_n<\lambda\}$. Then we have the following claims\\
(\romannumeral 1) if
$$\mathop{\lim\inf}\limits_{|x|\to \infty }|x|^{2s}c(x)\geq0,$$
there exists a constant $R_0>0$ (depending on c, but independent of $w$) such that
$$|x|<R_0.$$
(\romannumeral 2) if $c$ is bounded below in $\Sigma_\lambda$, there exists a constant $\ell>0$ (depending on the lower bound of $c$, but independent of $w$) such that
$$d(x,T_\lambda)> \ell,$$
where $T_\lambda=\{x\in\mathbb{R}^n|x_n=\lambda\}.$
\end{pro}
We want to point out that the constant $\ell$ is non-increasing about $\lambda$, since $\ell$ is non-decreasing about the lower bound of $c$,
which can be seen from the proof of Proposition \ref{L5} in \cite{DL}.

In order to apply the method of moving planes to prove the nonexistence, we need to establish the following estimate.
\begin{lem}\label{L4}
Let  $s\in(0,1)$ and $w(y)\in L_{2s}\cap C_{loc}^{1,1}(\mathbb{R}^n)$ satisfy $w(y)=-w(y^-_\lambda)$.  Assume there exists $x\in \Sigma_\lambda$ such that
$w(x)=\mathop{\inf}\limits_{\Sigma_\lambda }w(y)<0$. Then we have
$$(-\Delta)^sw(x)\leq C(n,s)\Big[w(x)d^{-2s}+\int_{\Sigma_\lambda}
(w(x)-w(y))\Big(\frac{1}{|x-y|^{n+2s}}-\frac{1}{|x-y^-_\lambda|^{n+2s}}\Big)dy\Big],$$
where $d=d(x,T_\lambda)$ and the constant $C(n,s)$ is positive.
\end{lem}
\proof Applying the definition of fractional Laplacian, we have
\begin{align}\label{020}
(-\Delta)^sw(x)&=C(n,s)\int_{\mathbb{R}^n}\frac{w(x)-w(y)}{|x-y|^{n+2s}}dy
\notag\\&=C(n,s)\int_{\Sigma_\lambda}\frac{w(x)-w(y)}{|x-y|^{n+2s}}dy+C(n,s)
\int_{\mathbb{R}^n\setminus\Sigma_\lambda}\frac{w(x)-w(y)}{|x-y|^{n+2s}}dy
\notag\\&=C(n,s)\int_{\Sigma_\lambda}\frac{w(x)-w(y)}{|x-y|^{n+2s}}dy+C(n,s)
\int_{\Sigma_\lambda}\frac{w(x)+w(y)}{|x-y^-_\lambda|^{n+2s}}dy
\notag\\&=C(n,s)\Big[\int_{\Sigma_\lambda}(w(x)-w(y))\Big(\frac{1}{|x-y|^{n+2s}}-\frac{1}
{|x-y^-_\lambda|^{n+2s}}\Big)dy+\int_{\Sigma_\lambda}\frac{2w(x)}{|x-y^-_\lambda|^{n+2s}}dy\Big].
\end{align}
By an elementary calculation (see \cite{CT}), we derive
$$
\int_{\Sigma_\lambda}\frac{2w(x)}{|x-y^-_\lambda|^{n+2s}}dy\cong C(n,s)w(x)d^{-2s}.
$$
Hence, combining this and (\ref{020}), we complete the proof of Lemma \ref{L4}.\endproof

In order to apply the method of moving spheres to prove the nonexistence, we need to establish a similar estimate as that of Lemma \ref{L4}. To this end, we need introduce some notations. For any real number $\lambda>0$, we denote
$$ S_\lambda=\{x\in\mathbb{R}^n \,|\, |x|=\lambda\},$$
$$ B_\lambda^+=B_\lambda^+(0)=\{|x|<\lambda\,|\,x_n>0\}.$$
Let $x^\lambda=\frac{\lambda^2x}{|x|^2}$ be the inversion of the point $ x=(x', x_n)$ about the sphere $S_\lambda$ and $x^*=(x', -x_n)$. Denote
$$ B_\lambda^-=\{x|x^*\in B_\lambda^+\},\ \ (B_\lambda^+)^C=\{x|x^\lambda\in B_\lambda^+\}, \ \ (B_\lambda^-)^C=\{x|x^\lambda\in B_\lambda^-\}.$$
\begin{lem}\label{L1}
Let $w(x)\in L_{2s}\cap C_{loc}^{1,1}(\mathbb{R}_+^n)$ satisfy
\begin{equation}\label{080}
w(x)=-w(x^*) \text{ and }  w(x)=-\big(\frac{\lambda}{|x|}\big)^{n-2s}w(x^\lambda), \ \ \forall x\in \mathbb{R}^n_+.
\end{equation}
Assume there exists $\tilde{x}\in B_\lambda^+$ such that $w(\tilde{x})=\inf_{B_\lambda^+}w(x)<0$. Then we have
$$
(-\Delta)^sw(\tilde{x})\leq C(n,s)
\Big[w(\tilde{x})\Big((\lambda-|\tilde{x}|)^{-2s}+\frac{\delta^n}{\tilde{x}_n^{n+2s}}\Big)+\int_{B_\lambda^+}
(w(\tilde{x})-w(y))h_\lambda(\tilde{x},y)dy\Big],$$
where $h_\lambda(\tilde{x},y)=\frac{1}{|\tilde{x}-y|^{n+2s}}-\frac{1}{|\frac{|y|}{\lambda}\tilde{x}-\frac{\lambda}{|y|}y|^{n+2s}}+
\frac{1}{|\frac{|y|}{\lambda}\tilde{x}-\frac{\lambda}{|y|}y^*|^{n+2s}}-\frac{1}{|\tilde{x}-y^*|^{n+2s}}>0$ for $\tilde{x},y\in B_\lambda^+$, $\delta=\min\{\tilde{x}_n,\lambda-|\tilde{x}|\}$ and $C(n,s)$ is a positive constant.
\end{lem}
\proof By the definition of fractional Laplacian and the assumptions (\ref{080}), we derive
\begin{align}\label{1}
(-\Delta)^sw(\tilde{x})&=C(n,s)\int_{\mathbb{R}^n}\frac{w(\tilde{x})-w(y)}{|\tilde{x}-y|^{n+2s}}dy
\notag\\&=C(n,s)\Big(\int_{B_\lambda^+}+\int_{(B_\lambda^+)^C}+\int_{B_\lambda^-}+\int_{(B_\lambda^-)^C}\Big)
\frac{w(\tilde{x})-w(y)}{|\tilde{x}-y|^{n+2s}}dy
\notag\\&=C(n,s)\Big(\int_{B_\lambda^+}\frac{w(\tilde{x})-w(y)}{|\tilde{x}-y|^{n+2s}}dy+
\int_{B_\lambda^+}\frac{(\frac{\lambda}{|y|})^{n-2s}w(\tilde{x})+w(y)}{|\frac{|y|}{\lambda}\tilde{x}-\frac{\lambda}{|y|}y|^{n+2s}}dy
\notag\\&+\int_{B_\lambda^-}\frac{w(\tilde{x})-w(y)}{|\tilde{x}-y|^{n+2s}}dy
+\int_{B_\lambda^-}\frac{(\frac{\lambda}{|y|})^{n-2s}w(\tilde{x})+w(y)}{|\frac{|y|}{\lambda}\tilde{x}-\frac{\lambda}{|y|}y|^{n+2s}}dy\Big)
\notag\\&=C(n,s)\Big[\int_{B_\lambda^+}(w(\tilde{x})-w(y))h_\lambda(\tilde{x},y)dy+\int_{B_\lambda^+}
\frac{(1+(\frac{\lambda}{|y|})^{n-2s})w(\tilde{x})}{|\frac{|y|}{\lambda}\tilde{x}-\frac{\lambda}{|y|}y|^{n+2s}}dy
\notag\\&+\int_{B_\lambda^+}\frac{2w(\tilde{x})}{|\tilde{x}-y^*|^{n+2s}}dy
+\int_{B_\lambda^+}\frac{(\frac{\lambda}{|y|})^{n-2s}-1)w(\tilde{x})}{|\frac{|y|}{\lambda}\tilde{x}-\frac{\lambda}{|y|}y^*|^{n+2s}}dy\Big].
\end{align}
Using similar arguments in \cite{LZ}, we can obtain that $h_\lambda(x,y)>0$ for $x,y\in B_\lambda^+$.

Furthermore, choose $r<\tilde{x}_n$ small such that $H:=\{x\in B_\delta(\tilde{x})|x_n>\tilde{x}_n\}\subset \{x\in\mathbb{R}^n|x_n>\tilde{x}_n\}\subset(B_r^+(0))^C$ where $\delta=\min\{\tilde{x}_n,\lambda-|\tilde{x}|\}$, then we calculate
\begin{align}\label{2}
\int_{B_\lambda^+}\frac{(1+(\frac{\lambda}{|y|})^{n-2s})w(\tilde{x})}{|\frac{|y|}{\lambda}\tilde{x}-\frac{\lambda}{|y|}y|^{n+2s}}dy&\leq
\int_{B_r^+}\frac{(1+(\frac{\lambda}{|y|})^{n-2s})w(\tilde{x})}{|\frac{|y|}{\lambda}\tilde{x}-\frac{\lambda}{|y|}y|^{n+2s}}dy
\notag\\&=\int_{(B_r^+)^C}\frac{(1+(\frac{\lambda}{|y^\lambda|})^{n-2s})w(\tilde{x})}{(\frac{|y^\lambda|}{\lambda})^{n+2s}|\tilde{x}-y|^{n+2s}}\big(\frac{\lambda}{|y|}\big)^{2n}dy
\notag\\&=w(\tilde{x})\int_{(B_r^+)^C}\frac{1}{|\tilde{x}-y|^{n+2s}}\big(1+(\frac{\lambda}{|y|})^{n-2s}\big)dy
\notag\\&\leq w(\tilde{x})\int_{\{x\in\mathbb{R}^n|x_n>\tilde{x}_n\}\setminus H}\frac{1}{|\tilde{x}-y|^{n+2s}}dy
\notag\\&\leq C(n) w(\tilde{x})\int_{\delta}^{+\infty}r^{-2s-1}dr
\notag\\&\leq C(n,s)w(\tilde{x})\delta^{-2s}
\notag\\&\leq C(n,s)w(\tilde{x})(\lambda-|\tilde{x}|)^{-2s}
\end{align}
From the definition $\delta=\min\{\tilde{x}_n,\lambda-|\tilde{x}|\}$, we have $|\tilde{x}-y^*|<C\tilde{x}_n$ for any $y\in B_\delta^+(\tilde{x})$. Simple calculations imply that
\begin{align}\label{3}
\int_{B_\lambda^+}\frac{2w(\tilde{x})}{|\tilde{x}-y^*|^{n+2s}}dy
\leq Cw(\tilde{x})\int_{B_\delta^+(\tilde{x})}\frac{1}{\tilde{x}_n^{n+2s}}dy
\leq C(n)w(\tilde{x})\frac{\delta^n}{\tilde{x}_n^{n+2s}}.
\end{align}
It is easy to see that
\begin{equation}\label{SS10}
\int_{B_\lambda^+}\frac{((\frac{\lambda}{|y|})^{n-2s}-1)w(\tilde{x})}{|\frac{|y|}{\lambda}\tilde{x}-\frac{\lambda}{|y|}y^*|^{n+2s}}\leq 0.
\end{equation}
Therefore, from (\ref{1})-(\ref{SS10}), we conclude the proof.\endproof

\begin{lem}\label {L2}
Let $\alpha,\beta >-n$. Suppose that $(u,v)$ is a nonnegative classical solution to the following system
\begin{equation}\label{4}
\left\{\begin{aligned}
&(-\Delta)^{s_1}u(x)=|x|^\alpha v^p(x),\ \,u(x)\geq 0,\ \ x\in\mathbb{R}_+^n,
\\&(-\Delta)^{s_2}v(x)=|x|^\beta u^q(x),\ \,v(x)\geq 0,\ \ x\in\mathbb{R}_+^n,
\\&u(x)=-u(x^*), \ \ v(x)=-v(x^*),\ \ \ \ x\in\mathbb{R}^n.
\end{aligned}\right.
\end{equation}
Then for $x\in \mathbb{R}^n_+$, we have
$$\left\{\begin{aligned}
&u(x)\geq C\int_{\mathbb{R}_+^n}\Big(\frac{1}{|x-y|^{n-2s_1}}-\frac{1}{|x^*-y|^{n-2s_1}}\Big)|y|^\alpha v^p(y)dy,
\\&v(x)\geq C\int_{\mathbb{R}_+^n}\Big(\frac{1}{|x-y|^{n-2s_2}}-\frac{1}{|x^*-y|^{n-2s_2}}\Big)|y|^\beta u^p(y)dy,
\end{aligned}\right.$$
where $C$ is a positive constant.
\end{lem}
{\bf Proof }
Define a cut off function $\eta(x)\in C_0^{\infty}(\mathbb{R}^n)$ satisfying $\eta(x)=0$ for $|x|>1$ and $\eta(x)=1$ for $|x|<\frac{1}{2}$. Denote $\eta_R(x)=\eta(\frac{x}{R})$ for large $R$ and
$$u_R(x)=C\int_{\mathbb{R}_+^n}\Big(\frac{1}{|x-y|^{n-2s_1}}-\frac{1}{|x^*-y|^{n-2s_1}}\Big)\eta_R(y)|y|^\alpha v^p(y)dy,$$
$$v_R(x)=C\int_{\mathbb{R}_+^n}\Big(\frac{1}{|x-y|^{n-2s_2}}-\frac{1}{|x^*-y|^{n-2s_2}}\Big)\eta_R(y)|y|^\beta u^q(y)dy.$$
Note that $(u_R(x), v_R(x))$ is a solution to the following system
\begin{equation}\label{5}
\left\{\begin{aligned}
&(-\Delta)^{s_1}u_R(x)=\eta_R(x)|x|^\alpha v^p(x),\ \ x\in\mathbb{R}_+^n,
\\&(-\Delta)^{s_2}v_R(x)=\eta_R(x)|x|^\beta u^q(x),\ \ x\in\mathbb{R}_+^n,
\\&u_R(x)=-u_R(x^*), \ \ v_R(x)=-v_R(x^*),\ \ \ \ x\in\mathbb{R}^n.
\end{aligned}\right.
\end{equation}
Let $U_R(x)=u(x)-u_R(x)$ and $V_R(x)=v(x)-v_R(x)$. From (\ref{4}) and (\ref{5}), we derive
\begin{equation}\label{b1}
\left\{\begin{aligned}
&(-\Delta)^{s_1}U_R(x)=|x|^\alpha v^p(x)-\eta_R(x)|x|^\alpha v^p(x)\geq 0,\ \ x\in\mathbb{R}_+^n,
\\&(-\Delta)^{s_2}V_R(x)=|x|^\beta u^q(x)-\eta_R(x)|x|^\beta u^q(x)\geq 0,\ \ x\in\mathbb{R}_+^n,
\\&U_R(x)=-U_R(x^*),\ \ V_R(x)=-V_R(x^*),\ \ \ \ x\in\mathbb{R}^n.
\end{aligned}\right.
\end{equation}
By the definitions of $U_R(x)$ and $V_R(x)$, obviously, for $x\in \mathbb{R}^n_+$,
\begin{equation}\label{070}
\mathop{\varliminf}\limits_{|x| \to \infty} U_R(x)\geq 0  \text { and } \mathop{\varliminf}\limits_{|x|\to \infty} V_R(x)\geq 0,
\end{equation}
where we used the assumptions $\alpha,\beta>-n$.

Next we claim that $U_R(x)\geq 0$ and $V_R(x)\geq 0 $ for $ x\in \mathbb{R}_+^n$. If not, from (\ref{070}) we know that there exists some $\hat{x}\in\mathbb{R}_+^n$ such that $U_R(\hat{x})=\mathop{\inf}\limits_{\mathbb{R}_+^n}U_R(x)<0$. Then

$$\begin{aligned}
(-\Delta)^{s_1}U_R(\hat{x})&=C(n,s_1)\int_{\mathbb{R}^n}\frac{U_R(\hat{x})-U_R(y)}{|\hat{x}-y|^{n+2s_1}}dy
\\&=C(n,s_1)\int_{\mathbb{R}_+^n}\frac{U_R(\hat{x})-U_R(y)}{|\hat{x}-y|^{n+2s_1}}dy+C(n,s_1)\int_{\mathbb{R}_+^n}\frac{U_R(\hat{x})
+U_R(y)}{|\hat{x}-y^*|^{n+2s_1}}dy
\\&=C(n,s_1)\Big[\int_{\mathbb{R}_+^n}(U_R(\hat{x})-U_R(y))\Big(\frac{1}{|\hat{x}-y|^{n+2s_1}}-\frac{1}{|\hat{x}-y^*|^{n+2s_1}}\Big)+\frac{2U_R(\hat{x})}
{|\hat{x}-y^*|^{n+2s_1}}\Big]dy
\\&<0.
\end{aligned}$$
This leads a contradiction with the first equation in (\ref{b1}). Thus, $U_R(x)\geq 0$ holds true for any $x\in \mathbb{R}_+^n$, that is, $u(x)\geq u_R(x)$ in $\mathbb{R}_+^n$. Letting $R \to \infty$, we obtain
$$u(x)\geq C\int_{\mathbb{R}_+^n}\Big(\frac{1}{|x-y|^{n-2s_1}}-\frac{1}{|x^*-y|^{n-2s_1}}\Big)|y|^\alpha v^p(y)dy.$$
Similarly, one has
$$v(x)\geq C\int_{\mathbb{R}_+^n}\Big(\frac{1}{|x-y|^{n-2s_2}}-\frac{1}{|x^*-y|^{n-2s_2}}\Big)|y|^\beta u^q(y)dy.$$

$\hfill\square$

\section{Proof of Theorem \ref{SST1}}
In this section, we are ready to prove  Theorem \ref{SST1}. To this end, for (i) of Theorem \ref{SST1}, namely the case $\min\{p+2s_1,p+2s_1+\alpha\}>1$ and $\{q+2s_2,q+2s_2+\beta\}>1$, we use the method of moving spheres to derive a lower bound for $u(x)$ and $v(x)$. Then, Lemma \ref{L2} and a ``Bootstrap" iteration process will give the better lower bounded estimates which can imply the nonexistence result. For  (ii) of Theorem \ref{SST1}, namely the case $0<pq<1$, $\alpha\geq -2s_1pq$ and $\beta\geq -2s_2pq$, a direct application of Lemma \ref{L2} and iteration technique may give its proof.

\subsection{ Proof of ($\romannumeral 1$) of Theorem \ref{SST1}}
\proof By contradiction, assume that $(u,v) \not\equiv (0,0)$, then we can derive that $u>0$ and $v>0$ in $\mathbb{R}^n_+$.
Indeed, if there exists some $\hat{x}\in\mathbb{R}_+^n$ such that $u(\hat{x})=0$, from the definition of $(-\Delta)^{s_1}$, we have

$$(-\Delta)^{s_1}u(\hat{x})=\int_{\mathbb{R}^n}\frac{-u(y)}{|\hat{x}-y|^{n+2s_1}}dy<0,$$
which contradicts with the equation
$$(-\Delta)^{s_1}u(\hat{x})=|\hat{x}|^\alpha v^p(\hat{x})\geq 0.$$
Thus $u(x)>0,$
and using the same arguments as above, we easily obtain $v(x)>0$. Therefore, we may assume that $u(x)>0$ and $v(x)>0$ in the rest proof of (i) of Theorem 1.

Let $u_\lambda(x)$ and $v_\lambda(x)$ be the Kelvin transform of $u(x)$ and $v(x)$ centered at origin respectively
$$u_\lambda(x)=\Big(\frac{\lambda}{|x|}\Big)^{n-2s_1}u\Big(\frac{\lambda^2x}{|x|^2}\Big),$$
$$v_\lambda(x)=\Big(\frac{\lambda}{|x|}\Big)^{n-2s_2}v\Big(\frac{\lambda^2x}{|x|^2}\Big)$$
for arbitrary $x\in\mathbb{R}^n\setminus \{0\}$.
By an elementary calculation, $u_\lambda(x)$ and $v_\lambda(x)$ satisfy the following system
$$\left\{\begin{aligned}
&(-\Delta)^{s_1}u_\lambda(x)=|x|^{\alpha}\Big(\frac{\lambda}{|x|}\Big)^{\tau_1}v_\lambda^p(x), \ \ x\in \mathbb{R}_+^n,
\\&(-\Delta)^{s_2}v_\lambda(x)=|x|^{\beta}\Big(\frac{\lambda}{|x|}\Big)^{\tau_2}u_\lambda^q(x), \ \ x\in \mathbb{R}_+^n,
\end{aligned}\right.
$$
where $$\tau_1=n+2s_1+2\alpha-p(n-2s_2)\ \ \text{and} \ \ \tau_2=n+2s_2+2\beta -q(n-2s_1).$$
Note that both $\tau_1$ and $\tau_2$ are nonnegative and they will not be zero simultaneously, since $(p,q)\in \mathcal{R}_{sub}$.

Denote
$$U_\lambda(x)=u_\lambda(x)-u(x)\ \ \text{and  }\ \ V_\lambda(x)=v_\lambda(x)-v(x).$$
By elementary calculations and the mean value theorem, for $x\in B_\lambda^+$, there holds
\begin{align}\label{S5}
(-\Delta)^{s_1}U_\lambda(x)&=|x|^{\alpha}\Big(\frac{\lambda}{|x|}\Big)^{\tau_1}v_\lambda^p(x)-|x|^{\alpha}v^p(x)
\\&=|x|^\alpha\Big[(v_\lambda^p(x)-v^p(x))+\Big(\Big(\frac{\lambda}{|x|}\Big)^{\tau_1}-1\Big)v_\lambda^p(x)\Big]\geq|x|^\alpha p\xi_{\lambda}^{p-1}(x)V_\lambda(x),
\notag\\(-\Delta)^{s_2}V_\lambda(x)&\geq|x|^\beta q\eta_{\lambda}^{q-1}(x)U_\lambda(x),
\end{align}
where $\xi_\lambda(x)$ is between $v(x)$ and $v_\lambda(x)$, $\eta_\lambda(x)$ is between $u(x)$ and $u_\lambda(x)$. Note that
\begin{equation}
U_\lambda(x)=-\Big(\frac{\lambda}{|x|}\Big)^{n-2s_1}U_\lambda(x^\lambda)\ \ \text{and}\ \ V_\lambda(x)=-\Big(\frac{\lambda}{|x|}\Big)^{n-2s_2}V_\lambda(x^\lambda).
\end{equation}

Next, we will use the method of moving spheres to claim that $U_\lambda(x)\geq 0$ and $V_\lambda(x)\geq 0$ in $B_\lambda^+$ for any $\lambda>0$.

\textbf{Step 1.} Give a start point. We show that for sufficiently small $\lambda>0$,
\begin{equation}\label{S7}
U_\lambda(x)\geq 0\ \ \text{and}\ \ V_\lambda(x)\geq 0, \ \ x\in B_\lambda^+.
\end{equation}
Suppose (\ref{S7}) is not true, there must exist a point $\bar{x}\in B_\lambda^+$ such that at least one of $U_\lambda(\bar{x})$ and $V_\lambda(\bar{x})$ is negative at this point. Without loss of generality, we assume
$$U_\lambda(\bar{x})=\inf_{x\in B_\lambda^+}\{U_\lambda(x),V_\lambda(x)\}<0.$$

We will obtain contradictions for all four possible cases respectively.

\textbf{Case 1.}  $(p,q)\in \mathcal{R}_{sub}$ and $ p\geq1, q\geq1$. Due to $p\geq 1$, then we can take $\xi_\lambda(x)=v(x)$ in (\ref{S5}). From the equation (\ref{S5}) and Lemma \ref{L1}, we have

\begin{equation}\label {S11}
\begin{aligned}
|\bar{x}|^\alpha pv^{p-1}(\bar{x})U_\lambda(\bar{x})&\leq|\bar{x}|^\alpha pv^{p-1}(\bar{x})V_\lambda(\bar{x})\leq(-\Delta)^{s_1}U_\lambda(\bar{x}) \\&\leq C(n,s_1)U_\lambda(\bar{x})\Big((\lambda-|\bar{x}|)^{-2s_1}+\frac{\delta^n}{\bar{x}_n^{n+2s_1}}\Big).
\end{aligned}\end{equation}
Hence,
\begin{equation}\label{001}
v^{p-1}(\bar{x})\geq\frac{C(n,s_1)\Big((\lambda-|\bar{x}|)^{-2s_1}+\frac{\delta^n}{\bar{x}_n^{n+2s_1}}\Big)}{p|\bar{x}|^\alpha }.
\end{equation}

If $\delta=\min\{\lambda-|\bar{x}|,\bar{x}_n\}=\lambda-|\bar{x}|$, which implies that $\lambda-|\bar{x}|\leq\bar{x}_n\leq|\bar{x}|$,
using the fact and (\ref{001}), we obtain
\begin{equation}\label{002}
v^{p-1}(\bar{x})\geq\frac{C(\lambda-|\bar{x}|)^{-2s_1}}{p|\bar{x}|^\alpha }\geq C|\bar{x}|^{-2s_1-\alpha}.
\end{equation}
As $\lambda\to 0$, the right of (\ref{002}) will go to infinity since $\alpha>-2s_1$. This is impossible.

If $\delta=\min\{\lambda-|\bar{x}|,\bar{x}_n\}=\bar{x}_n$, from $\bar{x}_n\leq|\bar{x}|$ and (\ref{001}), we derive
\begin{equation}\label{003}
v^{p-1}(\bar{x})\geq\frac{C\bar{x}_n^{-2s_1}}{p|\bar{x}|^\alpha }\geq C|\bar{x}|^{-2s_1-\alpha},
\end{equation}
which is also impossible.

\textbf{Case 2.} $0<p,q<1$. Due to $p<1$, we can take $\xi_\lambda(x) =v_\lambda(x)$. From the equation (\ref{S5}) and Lemma \ref{L1}, we have
\begin{equation}\label {S12}
\begin{aligned}
|\bar{x}|^\alpha p v_\lambda^{p-1}(\bar{x})U_\lambda(\bar{x})&\leq|\bar{x}|^\alpha p v_\lambda^{p-1}(\bar{x})V_\lambda(\bar{x})\leq (-\Delta)^{s_1}U_\lambda(\bar{x})
\\&\leq C(n,s_1)U_\lambda(\bar{x})\Big((\lambda-|\bar{x}|)^{-2s_1}+\frac{\delta^n}{\bar{x}_n^{n+2s_1}}\Big).
\end{aligned}
\end{equation}
Analogous to (\ref{002}) and (\ref{003}), there hold
\begin{equation}\label{004}
v_\lambda^{p-1}(\bar{x})\geq\frac{C\bar{x}_n^{-2s_1}}{|\bar{x}|^\alpha },
\end{equation}
or
\begin{equation}\label{014}
v_\lambda^{p-1}(\bar{x})\geq\frac{C(\lambda-|\bar{x}|)^{-2s_1}}{|\bar{x}|^\alpha }.
\end{equation}
Applying Lemma \ref{L2} and the mean value theorem, we obtain that for $x\in B_1^+$,
$$\begin{aligned}
u(x)&\geq C\int_{\mathbb{R}_+^n}\Big(\frac{1}{|x-y|^{n-2s_1}}-\frac{1}{|x^*-y|^{n-2s_1}}\Big)|y|^{\alpha}v^p(y)dy
\\&\geq C\int_{B_{1}(2e_n)}\Big(\frac{1}{|x-y|^{n-2s_1}}-\frac{1}{|x^*-y|^{n-2s_1}}\Big)dy
\\&\geq C\int_{B_{1}(2e_n)}\frac{x_ny_n}{|x^*-y|^{n-2s_1+2}}dy
\\&\geq Cx_n.
\end{aligned}$$
For $x\in (B_1^+)^C$, we derive
$$\begin{aligned}
u(x)&\geq C\int_{\mathbb{R}_+^n}\Big(\frac{1}{|x-y|^{n-2s_1}}-\frac{1}{|x^*-y|^{n-2s_1}}\Big)|y|^{\alpha}v^p(y)dy
\\&\geq C\int_{B_{1}(2e_n)}\frac{x_ny_n}{|x^*-y|^{n-2s_1+2}}dy
\\&\geq C\frac{x_n}{|x|^{n-2s_1+2}}.
\end{aligned}$$
Similarly, we have
$$v(x)\geq\left\{\begin{aligned}
 &Cx_n, &x\in B_1^+,
\\&C\frac{x_n}{|x|^{n-2s_2+2}},&x\in (B_1^+)^C.
\end{aligned}\right.$$
Then by the definition of $v_\lambda(x)$, we obtain
\begin{equation}\label{015}
v_\lambda(x)\geq \left\{\begin{aligned}
 &C\Big(\frac{\lambda}{|x|}\Big)^{n-2s_2+2}x_n,&x^\lambda\in B_1^+,
\\&C\frac{x_n}{\lambda^{n-2s_2+2}},&x^\lambda\in (B_1^+)^C.
\end{aligned}\right.
\end{equation}

If $\delta=\lambda-|\bar{x}|\leq\bar{x}_n$, that is, $|\bar{x}|\geq\frac{\lambda}{2}$, combining (\ref{014}) and (\ref{015}),we conclude that  for $\bar{x}\in B_\lambda^+$ and sufficiently small $\lambda$,
$$(\lambda-|\bar{x}|)^{-2s_1} \leq C\bar{x}_n^{p-1}|\bar{x}|^\alpha\leq C(\lambda-|\bar{x}|)^{p-1}|\bar{x}|^\alpha,$$
which gives that
\begin{equation}\label{007}
\Big(\frac{\lambda}{|\bar{x}|}-1\Big)^{-2s_1-p+1}\leq C|\bar{x}|^{p+2s_1+\alpha-1}.
\end{equation}
Due to $\min\{p+2s_1, p+2s_1+\alpha\}>1$ and $\frac{\lambda}{|\bar{x}|}\leq2$, the inequality (\ref{007}) is impossible as $\lambda>0$ sufficiently small.

If $\delta=\bar{x}_n$, it follows (\ref{004}) and (\ref{015}) that
$$\frac{C\bar{x}_n^{-2s_1}}{|\bar{x}|^\alpha }\leq v_\lambda^{p-1}(\bar{x})\leq C \bar{x}_n^{p-1},$$
which implies that
$$|\bar{x}|^{-p-2s_1-\alpha+1}\leq C \ \ \text{ if } \alpha\geq 0,\ \ \text {and}\ \  \bar{x}_n^{-p-2s_1-\alpha+1}\leq C \ \ \text{ if } \alpha< 0.$$
These arguments yield a contradiction since the left terms go to infinity as $\min\{p+2s_1,p+2s_1+\alpha\}>1$ and $\lambda>0$ small enough.

For Case 3: $(p,q)\in \mathcal{R}_{sub},  p\geq1, 0<q<1$ and Case 4: $(p,q)\in \mathcal{R}_{sub}, 0<p<1,  q\geq1$, similar argument as that of Case 1 and Case 2 can show that $U_\lambda(x)\geq0$ and $V_\lambda(x)\geq 0$ in $B_\lambda^+$ for sufficiently small $\lambda>0$. Therefore, (\ref{S7}) holds.

\textbf{Step 2.} Now we move the sphere $S_\lambda$ outward as long as (\ref{S7}) holds. Define
$$\lambda_0=\sup\{\lambda|\,U_\mu(x)\geq 0, V_\mu(x)\geq 0,\ \  x\in B_\mu^+, \ \ \forall\ \ 0<\mu<\lambda\}.$$
We will show that $\lambda_0=+\infty$. Suppose on the contrary that $0<\lambda_0<+\infty$. We want to show that there exists some small $\varepsilon>0$ such that for any $\lambda\in (\lambda_0,\lambda_0+\varepsilon)$,
\begin{equation}\label{S15}
U_\lambda(x)\geq 0\ \ \text{and}\ \ V_\lambda(x)\geq 0, \ \ x\in B_\lambda^+.
\end{equation}
This implies that the plane $S_{\lambda_0}$ will be moved outward a little bit further, which contradicts with the definition of $\lambda_0$.

Firstly, we claim that
\begin{equation}\label{S16}
U_{\lambda_0}(x)> 0\ \ \text{and}\ \ V_{\lambda_0}(x)> 0, \ \ x\in B_{\lambda_0}^+.
\end{equation}
Indeed, if there exists some point $x^0\in B_{\lambda_0}^+$ such that $U_{\lambda_0}(x^0)= 0$, we have
\begin{equation}\label{012}
(-\Delta)^{s_1}U_{\lambda_0}(x^0)=C\int_{\mathbb{R}^n}\frac{-U_{\lambda_0}(y)}{|x^0-y|^{n+2s_1}}dy\leq 0.
\end{equation}
On the other hand, it is easy to get that
\begin{equation}\label{013}
\begin{aligned}
&(-\Delta)^{s_1}U_{\lambda_0}(x^0)=|x^0|^{\alpha}\Big(\frac{\lambda_0}{|x^0|}\Big)^{\tau_1}v_{\lambda_0}^p(x^0)-|x^0|^{\alpha}v^p(x^0)
\\&=|x^0|^\alpha \Big(\Big(\frac{\lambda_0}{|x^0|}\Big)^{\tau_1}-1\Big)v_\lambda^p(x^0)+p|x^0|^\alpha \xi_\lambda^{p-1}(x^0)V_{\lambda_0}(x^0).
\end{aligned}
\end{equation}
If $\tau_1>0$, then $(-\Delta)^{s_1}U_{\lambda_0}(x^0)>0$, where we use the fact $V_{\lambda_0}\geq 0$ and .
If $\tau_1=0$, then we have that $\tau_2>0$ and $V_{\lambda_0}(x^0)=0$ follows from (\ref{013}). Using an argument similar to (\ref{012}) and (\ref{013}), we derive
$$\begin{aligned}0\geq(-\Delta)^{s_2}V_{\lambda_0}(x^0)&= |x^0|^\beta \Big(\Big(\frac{\lambda_0}{|x^0|}\Big)^{\tau_2}-1\Big)u_\lambda^q(x^0)+|x^0|^\beta \eta_\lambda^{q-1}(x^0)U_{\lambda_0}(x^0)
\\&= |x^0|^\beta \Big(\Big(\frac{\lambda_0}{|x^0|}\Big)^{\tau_2}-1\Big)u_\lambda^q(x^0)>0,\end{aligned}$$
which is absurd. Thus $U_{\lambda_0}(x)> 0$ is proved. Using the fact and similar arguments as above, we derive that $V_{\lambda_0}(x)> 0$. Hence, (\ref{S16}) holds.

Next, we will show that the sphere can be moved further outward. The continuity of $u(x)$ and (\ref{S16}) yield that there exists some sufficiently small $l\in(0,\frac{\lambda_0}{2})$ and $\varepsilon_1\in(0,\frac{\lambda_0}{2})$ such that for $\lambda\in(\lambda_0,\lambda_0+\varepsilon_1)$,
\begin{equation}\label{S17}
U_\lambda(x)\geq 0, \ \ \ \, x\in B_{\lambda_0-l}^+.
\end{equation}
For $x\in B_\lambda^+\setminus B_{\lambda_0-l}^+$, using the similar proof of (\ref{S7}), we can deduce that
\begin{equation}\label{S19}
U_\lambda(x)\geq 0, \ \ \ \ x\in B_\lambda^+\setminus B_{\lambda_0-l}^+ .
\end{equation}
Note that the distance between $\bar{x}$ and $S_\lambda$, i.e. $\lambda-|\bar{x}|$, plays an important role in this process.

Hence, It follows from (\ref{S17}) and (\ref{S19}) that for all $\lambda\in (\lambda_0,\lambda_0+\varepsilon_1)$,  $$U_\lambda(x)\geq 0,\ \ x\in B_\lambda^+.$$
Similarly, we can also prove that for all $\lambda\in (\lambda_0,\lambda_0+\varepsilon_2)$,
$$V_\lambda(x)\geq 0,\ \ x\in B_\lambda^+.$$
Let $\varepsilon=\min\{\varepsilon_1,\varepsilon_2\}$, therefore, (\ref{S15}) can be completely concluded. This contradicts with the definition of $\lambda_0$. So $\lambda_0=+\infty$.

Then, we have for every $\lambda>0$,
$$U_\lambda(x)\geq 0,\ \ V_\lambda(x)\geq 0, \ \  x\in B_\lambda^+,$$
which gives that,
\begin{equation}\label{SS15}
u(x)\geq\Big(\frac{\lambda}{|x|}\Big)^{n-2s_1}u\Big(\frac{\lambda^2x}{|x|^2}\Big), \ \ \ \ \forall \,|x|\geq\lambda, \ \ x\in\mathbb{R}_+^n, \  \ \forall\, 0<\lambda<+\infty.
\end{equation}
\begin{equation}\label{SS16}
v(x)\geq\Big(\frac{\lambda}{|x|}\Big)^{n-2s_2}v\Big(\frac{\lambda^2x}{|x|^2}\Big), \ \ \ \ \forall \,|x|\geq\lambda, \ \ x\in\mathbb{R}_+^n, \  \ \forall\, 0<\lambda<+\infty.
\end{equation}
For any given $|x|\geq 1$, let $\lambda=\sqrt{|x|}$, then it follows from (\ref{SS15}) that
\begin{equation}\label{16}
u(x)\geq\Big(\min_{x\in S_1^+}u(x)\Big)\frac{1}{|x|^{\frac{n-2s_1}{2}}}:=\frac{C}{|x|^{\frac{n-2s_1}{2}}}\geq\frac{Cx_n}{|x|^{\frac{n-2s_1}{2}+1}},
\end{equation}
and similarly from (\ref{SS16}), we obtain
\begin{equation}\label{16}
v(x)\geq\frac{Cx_n}{|x|^{\frac{n-2s_2}{2}+1}}.
\end{equation}

Now we make full use of the above properties to derive some lower bound estimates of solutions to (\ref{SI3}) through iteration technique.

Let $\theta_0=\frac{n-2s_1}{2}+1$, $\sigma_0=\frac{n-2s_2}{2}+1$. From Lemma \ref{L2}, inequality (\ref{16}) and the mean value theorem, we have for $x_n>1$,
\begin{align}\label{17}
u(x)&\geq C\int_{\mathbb{R}_+^n}\Big(\frac{1}{|x-y|^{n-2s_1}}-\frac{1}{|x^*-y|^{n-2s_1}}\Big)|y|^{\alpha}v^p(y)dy
\notag\\&\geq C\int_{2|x|}^{4|x|}\int_{2|x|\leq |y'|\leq 4|x|}\frac{x_ny_n}{|x^*-y|^{n-2s_1+2}}\frac{y_n^p}{|y|^{p\sigma_0-\alpha}}dy'dy_n
\notag\\&\geq C\frac{x_n}{|x|^{(n-2s_1+2)+(p\sigma_0-\alpha)}}\int_{2|x|}^{4|x|}\int_{2|x|\leq |y'|\leq 4|x|}y_n^{p+1}dy'dy_n
\notag\\&\geq C\frac{x_n}{|x|^{p\sigma_0-\alpha-2s_1+3-(p+2)}}.
\end{align}
Similarly, we have
\begin{equation}\label{SS18}
v(x) \geq C\frac{x_n}{|x|^{q(\theta_0-1)-(\beta+2s_2)+1}}
\end{equation}
Denote $\theta_1=p(\sigma_0-1)-(\alpha+2 s_1)+1$ and $\sigma_1=q(\theta_0-1)-(\beta+2s_2)+1$. Repeat the above process replacing (\ref{16}) by (\ref{SS18}),  then we have
$$u(x)\geq  C\int_{2|x|}^{4|x|}\int_{2|x|\leq |y'|\leq 4|x| }\frac{x_ny_n}{|x^*-y|^{n-2s_1+2}}\frac{y_n^p}{|y|^{p\sigma_1-\alpha}}dy'dy_n\geq C\frac{x_n}{|x|^{\theta_2}},$$
and analogously,
\begin{equation}\label{S44}
v(x)\geq C\frac{x_n}{|x|^{\sigma_2}},
\end{equation}
where $\theta_2=p(\sigma_1-1)-\alpha-2s_1+1$ and $\sigma_2=q(\theta_1-1)-\beta-2s_2+1$.

After such $k$ iteration steps, we derive
\begin{equation}\label {SS33}
u(x)\geq  C\frac{x_n}{|x|^{\theta_{k+1}}}, \ \ \ \ v(x)\geq C\frac{x_n}{|x|^{\sigma_{k+1}}},
\end{equation}
where $\theta_{k+1}=p(\sigma_{k}-1)-\alpha-2s_1+1$ and $\sigma_{k+1}=q(\theta_{k}-1)-\beta-2s_2+1$. Elementary calculations give that

\begin{equation}\label{SS45}
\begin{aligned}&\theta_{2m}=\frac{n-2s_1}{2}(pq)^{m}-\Big[\big(p(2s_2+\beta)+(2s_1+\alpha)\big)\frac{1-(pq)^m}{1-pq}\Big]+1,
\\&\theta_{2m+1}=\Big(\frac{p(n-2s_2)}{2}-2s_1-\alpha \Big)(pq)^m-\Big[\big(p(2s_2+\beta)+(2s_1+\alpha)\big)\frac{1-(pq)^m}{1-pq}\Big]+1,
\\&\sigma_{2m}=\frac{n-2s_2}{2}(pq)^{m}-\Big[\big(q(2s_1+\alpha)+(2s_2+\beta)\big)\frac{1-(pq)^m}{1-pq}\Big]+1,
\\&\sigma_{2m+1}=\Big(\frac{q(n-2s_1)}{2}-2s_2-\beta \Big)(pq)^m-\Big[\big(q(2s_1+\alpha)+(2s_2+\beta)\big)\frac{1-(pq)^m}{1-pq}\Big]+1,
\end{aligned}
\end{equation}
where $m=0,1,2,\cdots$.

For the case $pq\geq 1$, we claim that both $\{\theta_k\}$ and $\{\sigma_k\}$ are decreasing sequences and unbounded from below. Denote
$$A_e=p\frac{n-2s_2}{2}-\frac{n-2s_1}{2}-2s_1-\alpha,$$
$$A_o=q\frac{n-2s_1}{2}-\frac{n-2s_2}{2}-2s_2-\beta.$$
Note that $A_e,A_o\leq 0$ and they will not be zero simultaneously, since $(p,q)\in \mathcal{R}_{sub}$.
By elementary calculations, we have
\begin{equation}\label{088}
\theta_{k+1}-\theta_{k}=\left\{\begin{aligned}
&p^{[\frac{k+1}{2}]}q^{[\frac{k}{2}]}A_e\leq 0,~~~  k \mbox{ is even},\\&
p^{[\frac{k+1}{2}]}q^{[\frac{k}{2}]}A_o\leq 0,~~~  k \mbox{ is odd}.
\end{aligned}\right.
\end{equation}
Hence the sequence $\{\theta_k\}$ is decreasing.

Combining the above properties of $A_e, A_o$, the assumption $pq\geq1$  and (\ref{088}), we deduce that $\theta_k\to -\infty$  as $k\to+\infty$. Similarly, we can show that $\{\sigma_k\}$ is also decreasing and $\sigma_k\to-\infty$. These and (\ref{SS33}) indicate that $u(x)$ and $v(x)$ don't belong to any $L_{2s}$. Hence, the solutions $(u,v)\equiv (0,0)$  when $pq\geq1$.

Next we consider the case $pq<1$. From (\ref{SS45}), we can conclude that
\begin{equation}\label {S47}
\theta_k\to1-\frac{p(2s_2+\beta)+(2s_1+\alpha)}{1-pq}, \ \ \ \ \sigma_k\to1-\frac{q(2s_1+\alpha)+(2s_2+\beta)}{1-pq}\ \ \ \ \text{as } k\to \infty.
\end{equation}
Since $p+2s_1+\alpha>1$ and $q+2s_2+\beta>1$, we have
$$\frac{p(2s_2+\beta)+(2s_1+\alpha)}{1-pq}-1>0, \ \  \ \ \frac{q(2s_1+\alpha)+(2s_2+\beta)}{1-pq}-1>0.$$
Hence, from (\ref{S47}) and (\ref{SS33}), we have for $x_n>1$,
\begin{equation}\label{S48}
u(x)\geq cx_n^{\frac{p(2s_2+\beta)+(2s_1+\alpha)}{1-pq}+o(1)}, \ \ \ \ v(x)\geq cx_n^{\frac{q(2s_1+\alpha)+(2s_2+\beta)}{1-pq}+o(1)}.
\end{equation}
Combining this with Lemma \ref{L2}, we have
$$\begin{aligned}
u(x)&\geq C\int_{\mathbb{R}_+^n}\Big(\frac{1}{|x-y|^{n-2s_1}}-\frac{1}{|x^*-y|^{n-2s_1}}\Big)|y|^{\alpha}v^p(y)dy
\\&\geq C \int_{2x_n}^{+\infty}\int_{\mathbb{R}^{n-1}}\Big(\frac{1}{|x-y|^{n-2s_1}}-\frac{1}{|x^*-y|^{n-2s_1}}\Big)
y_n^{(\frac{q(2s_1+\alpha)+(2s_2+\beta)}{1-pq}+o(1))p+\alpha}dy'dy_n
\\&\geq C\int_{2x_n}^{+\infty}y_n^{(\frac{q(2s_1+\alpha)+(2s_2+\beta)}{1-pq}+o(1))p+\alpha}dy_n\int_{\mathbb{R}^{n-1}}
\frac{x_ny_n}{\big(|x'-y'|^2+|x_n+y_n|^2\big)^{\frac{n-2s_1+2}{2}}}
dy'.
\end{aligned}$$
Let $x'-y'=(x_n+y_n)z'$, we have
$$\begin{aligned}
u(x)&\geq C\int_{2x_n}^{+\infty}\frac{y_n^{(\frac{q(2s_1+\alpha)+(2s_2+\beta)}{1-pq}+o(1))p+\alpha+1}x_n}{(x_n+y_n)^{3-2s_1}}dy_n
\int_{\mathbb{R}^{n-1}}\frac{1}{\big(|z'|^2+1\big)^{\frac{n-2s_1+2}{2}}}
dz'
\\&\geq C\int_{2x_n}^{+\infty}\frac{y_n^{(\frac{q(2s_1+\alpha)+(2s_2+\beta)}{1-pq}+o(1))p+\alpha+1}x_n}{(x_n+y_n)^{3-2s_1}}dy_n.
\end{aligned}$$
Let $y_n=x_nz_n$, one has
\begin{align}\label{SS38}
u(x)&\geq Cx_n^{(\frac{q(2s_1+\alpha)+(2s_2+\beta)}{1-pq}+o(1))p+2s_1+\alpha}\int_{2}^{+\infty}\frac{z_n^{\frac{q(2s_1+\alpha)
+(2s_2+\beta)}{1-pq}p+\alpha+1}}{(1+z_n)^{3-2s_1}}dz_n
\notag\\&\cong\int_{2}^{+\infty}\frac{1}{(z_n)^{3-2s_1-(\frac{q(2s_1+\alpha)+(2s_2+\beta)}{1-pq}p+\alpha+1)}}dy_n.
\end{align}
Due to $p+2s_1+\alpha>1$ and $q+2s_2+\beta>1$, we have $\frac{q(2s_1+\alpha)+(2s_2+\beta)}{1-pq}>1$, which implies that $$3-2s_1-\left(\frac{q(2s_1+\alpha)+(2s_2+\beta)}{1-pq}p+\alpha+1\right)<1.$$
This yields that the right hand side of (\ref{SS38}) is infinity, which is impossible. Hence, for $pq<1$ we also obtain $(u,v)\equiv (0,0).$

In sum, we conclude that there is no nontrivial classical solutions to system (\ref{SI3})  for the case that $(p,q)\in \mathcal{R}_{sub}$, $\min\{p+2s_1, p+2s_1+\alpha\}>1$, $\min\{q+2s_2,q+2s_2+\beta\}>1$.
\endproof

\subsection{Proof of (\romannumeral 2) of Theorem \ref{SST1}}
\begin{proof}Assume that $(u,v) \not\equiv (0,0)$, from the proof of (\romannumeral 1), we know that $u>0$ and $v>0$ in $\mathbb{R}_+^n$. Applying Lemma \ref{L2},  for $x_n>\min\{2,\frac{|x'|}{10}\}$, we have
\begin{equation}\label{SS39}
\begin{aligned}
u(x)&\geq C\int_{\mathbb{R}_+^n}\Big(\frac{1}{|x-y|^{n-2s_1}}-\frac{1}{|x^*-y|^{n-2s_1}}\Big)|y|^{\alpha}v^p(y)dy
\\&\geq C\int_{B_1(2e_n)}\Big(\frac{x_ny_n}{|x^*-y|^{n-2s_1+2}}\Big)|y|^{\alpha}v^p(y)dy
\\&\geq C\frac{1+|x|}{(1+|x|)^{n-2s_1+2}}\int_{B_1(2e_n)}y_n|y|^{\alpha}v^p(y)dy
\\&\geq C\frac{1}{(1+|x|)^{n-2s_1+1}}.
\end{aligned}
\end{equation}
Similarly, for $x_n>\min\{2,\frac{|x'|}{10}\}$, we have
\begin{equation}\label{S52}
v(x)\geq \frac{C}{(1+|x|)^{n-2s_2+1}}.
\end{equation}
Iterating  (\ref{SS39}) with (\ref{S52}), for $x_n>\min\{2,\frac{|x'|}{10}\}$, we get

\begin{align}\label{S53}
u(x)&\geq C\int_{B_{|x|}(0,4|x|)}\Big(\frac{1}{|x-y|^{n-2s_1}}-\frac{1}{|x^*-y|^{n-2s_1}}\Big)\frac{|y|^\alpha}{(1+|y|)^{p(n-2s_2+1)}}dy
\notag\\&\geq C\int_{B_{3|x|}(0,4|x|)\setminus B_{2|x|}(0,4|x|)}\frac{x_ny_n}{|x^*-y|^{n-2s_1+2}}\frac{|y|^\alpha}{(1+|y|)^{p(n-2s_2+1)}}dy
\notag\\&\geq C\frac{(1+|x|)^{2+\alpha}}{(1+|x|)^{n-2s_1+2+(n-2s_2+1)p}}\int_{B_{3|x|}(0,4|x|)\setminus B_{2|x|}(0,4|x|)}dy
\notag\\&\geq C\frac{1}{(1+|x|)^{(n-2s_2+1)p-2s_1-\alpha}},
\end{align}
Using the same argument as that of (\ref{S53}), for $x_n>\min\{2,\frac{|x'|}{10}\}$, we obtain
$$v(x)\geq C\frac{1}{(1+|x|)^{(n-2s_1+1)q-2s_2-\beta}}.$$
After $k$ iteration steps, it is easy to see that for $|x|$ large and $x_n>\min\{2,\frac{|x'|}{10}\}$,
$$u(x)\geq \frac{C}{(1+|x|)^{\gamma_k}}, \ \ \ \ v(x)\geq \frac{C}{(1+|x|)^{\delta_k}}.$$
Here
$$\gamma_k=\delta_{k-1}p-2s_1-\alpha,\ \ \ \ \delta_k=\gamma_{k-1}q-2s_2-\beta,$$
where
$$\gamma_1=(n-2s_2+1)p-2s_1-\alpha,\ \ \ \ \delta_1=(n-2s_1+1)q-2s_2-\beta.$$
Simple calculations imply that
$$\begin{aligned}
&\gamma_{2m}=(n-2s_1+1)(pq)^m-\Big[\big(p(2s_2+\beta)+(2s_1+\alpha)\big)\frac{1-(pq)^m}{1-pq}\Big],
\\&\gamma_{2m+1}=\big[p(n-2s_2+1)-(2s_1+\alpha)\big](pq)^m-\Big[\big(p(2s_2+\beta)+(2s_1+\alpha)\big)\frac{1-(pq)^m}{1-pq}\Big],
\\&\delta_{2m}=(n-2s_2+1)(pq)^m-\Big[\big(q(2s_1+\alpha)+(2s_2+\beta)\big)\frac{1-(pq)^m}{1-pq}\Big],
\\&\delta_{2m}=\big[q(n-2s_1+1)-(2s_2+\beta)\big](pq)^m-\Big[\big(q(2s_1+\alpha)+(2s_2+\beta)\big)\frac{1-(pq)^m}{1-pq}\Big],
\end{aligned}$$
where $m=0,1,2,\cdots.$

From $0<pq<1$, we obtain
$$\gamma_k\to-\frac{p(2s_2+\beta)+(2s_1+\alpha)}{1-pq},\ \ \ \ \delta_k\to-\frac{q(2s_1+\alpha)+(2s_2+\beta)}{1-pq},\ \ \ \ \text{as } k\to\infty.$$
This yields that for $|x|$ large,
$$u(x)\geq C(1+|x|)^{\frac{p(2s_2+\beta)+(2s_1+\alpha)}{1-pq}-o(1)}, \ \ \ \ v(x)\geq C(1+|x|)^{\frac{q(2s_1+\alpha)+(2s_2+\beta)}{1-pq}-o(1)},$$
as $k\to+\infty$. Due to $0<pq<1$ and $\alpha\geq-2s_1pq$, $\beta\geq-2s_2pq$, we have
$$\frac{p(2s_2+\beta)+(2s_1+\alpha)}{1-pq}>2s_1,\ \ \ \ \frac{q(2s_1+\alpha)+(2s_2+\beta)}{1-pq}>2s_2$$
This contradicts with the assumptions that $u(x)\in L_{2s_1}$ and $v(x)\in L_{2s_2}$. Hence, $(u,v)\equiv (0,0)$.\end{proof}

\section{Proof of Theorem \ref{T2}}

In this section, we give the proof of Theorem \ref{T2} by using the method of moving planes.
\proof Suppose on the contrary that $(u,v) \not\equiv (0,0)$, we know that $u>0$ and $v>0$ in $\mathbb{R}_+^n$.
Let $\bar{u}(x)$ and $\bar{v}(x)$ be the Kelvin transform of $u(x)$ and $v(x)$ centered at origin respectively
$$\bar{u}(x)=\Big(\frac{1}{|x|}\Big)^{n-2s_1}u\Big(\frac{x}{|x|^2}\Big),$$
$$\bar{v}(x)=\Big(\frac{1}{|x|}\Big)^{n-2s_2}v\Big(\frac{x}{|x|^2}\Big)$$
for arbitrary $x\in\mathbb{R}^n\setminus \{0\}$. Then, $\bar{u}(x)$ and $\bar{v}(x)$ satisfy the following system
$$\left\{\begin{aligned}
&(-\Delta)^{s_1}\bar{u}(x)=\Big(\frac{1}{|x|}\Big)^{\bar{\tau}_1}\bar{v}^p(x), \ \ x\in \mathbb{R}_+^n,
\\&(-\Delta)^{s_2}\bar{v}(x)=\Big(\frac{1}{|x|}\Big)^{\bar{\tau}_2}\bar{u}^q(x), \ \ x\in \mathbb{R}_+^n,
\\&\bar{u}(x)> 0, \ \ \bar{v}(x)> 0, \ \  x\in \mathbb{R}_+^n,
\\&\bar{u}(x',x_n)=-\bar{u}(x',-x_n),\ \ \bar{v}(x',x_n)=-\bar{v}(x',-x_n),\ \ \ \ x=(x',x_n)\in\mathbb{R}^n,
\end{aligned}\right.
$$
where
$$\bar{\tau}_1=n+2s_1+\alpha-p(n-2s_2)\ \ \text{ and }\ \ \bar{\tau}_2=n+2s_2+\beta -q(n-2s_1).$$
Note that $\bar{\tau}_1\leq0$ and $\bar{\tau}_2\leq 0$ due to $p\geq \frac{n+2s_1+\alpha}{n-2s_2}$ and $q\geq \frac{n+2s_2+\beta}{n-2s_1}$. Obviously, for $|x|$ large enough,
\begin{equation}\label{030}
\bar{u}(x)=O\Big(\frac{1}{|x|^{n-2s_1}}\Big)\ \ \text { and } \ \ \bar{v}(x)=O\Big(\frac{1}{|x|^{n-2s_2}}\Big).
\end{equation}

For any real number $\rho>0$ and $x\in \Sigma_\rho$,
define
$$\bar{U}_\rho(x)=\bar{u}(x^-_\rho)-\bar{u}(x), \ \  \bar{V}_\rho(x)=\bar{v}(x^-_\rho)-\bar{v}(x),$$
where $\Sigma_\rho$, $T_\rho$ and $x^-_\rho$ are defined as in Proposition \ref{L5}.
Then, due to $p,q>1$, we have
\begin{align}\label{031}
&(-\Delta)^{s_1}\bar{U}_\rho(x)=|x^-_\rho|^{-\bar{\tau}_1} |\bar{v}(x^-_\rho)|^{p-1}\bar v(x'_\rho)-|x|^{-\bar{\tau}_1} |\bar v(x)|^{p-1}\bar v(x)\geq |x|^{-\bar{\tau}_1} p|\bar{v}|^{p-1}(x)\bar{V}_\rho(x),
\\&(-\Delta)^{s_2}\bar{V}_\rho(x)\geq |x|^{-\bar{\tau}_2} q|\bar{U}_\rho|^{q-1}(x)\bar{U}_\rho(x).
\end{align}

\textbf{Step 1} We will claim that for $\rho>0$ sufficiently small,
\begin{equation}\label{032}
\bar{U}_\rho(x)\geq 0\ \ \text{and  }\ \ \bar{V}_\rho(x)\geq 0,\ \ \ \ x\in\Sigma_\rho.
\end{equation}
Otherwise, from (\ref{030}), there exists some $\bar{x}\in\Sigma_\rho\cap \mathbb{R}_+^n$ such that at leat one of $\bar{U}_\rho(x)$, $\bar{V}_\rho(x)$ is negative. Without loss of generality, we may assume that
$$\bar{U}_\rho(\bar{x})=\inf_{\Sigma_\rho}\{\bar{U}_\rho(x),\bar{V}_\rho(x)\}<0.$$

Combining the equation (\ref{031}) and Lemma \ref{L4}, we deduce
\begin{equation}\label{201}
\begin{aligned}
|\bar{x}|^{-\bar{\tau}_1} p\bar{v}^{p-1}(\bar{x})\bar{U}_\rho(\bar{x})&\leq|\bar{x}|^{-\bar{\tau}_1}
p\bar{v}^{p-1}(\bar{x})\bar{V}_\rho(\bar{x})\leq(-\Delta)^{s_1}\bar{U}_\rho(\bar{x})
\\& \leq C(n,s_1)\bar{U}_\rho(\bar{x})(\rho-\bar{x}_n)^{-2s_1}.
\end{aligned}
\end{equation}
This yields that
\begin{equation}\label{034}
|\bar{x}|^{-\bar{\tau}_1} p\bar{v}^{p-1}(\bar{x})\geq C(\rho-\bar{x}_n)^{-2s_1}\geq C\rho^{-2s_1}.
\end{equation}
Observe that (\ref{030}) and the decay conditions of $u$ and $v$ in Theorem \ref{T2} ensure that
$$
\mathop{\lim}\limits_{ x \to \infty}|x|^{-\bar{\tau}_1} p\bar{v}^{p-1}(x)=0, ~~\mathop{\overline{\lim}}\limits_{x\to 0}|x|^{-\bar{\tau}_1} p\bar{v}^{p-1}(x)\leq C,
$$
where we used the assumption $\alpha>-2s_2$. Hence, the inequality (\ref{034}) is impossible as $\rho>0$ is sufficiently small. Therefore, (\ref{032}) holds.

\textbf{Step 2} Move the plane $T_\rho$ upward along $x_n$-axis as long as (\ref{032}) holds. Let
$$\rho_0=\sup\{\rho \,|\, \bar{U}_\mu(x)\geq 0, \bar{V}_\mu(x)\geq 0, x\in \Sigma_\mu,\mu\leq \rho,\rho>0\}.$$
We will show that $\rho_0=+\infty$ by contradiction arguments.

Suppose on the contrary that $0<\rho_0<+\infty$. We will verify that
\begin{equation}\label{S26}
\bar{U}_{\rho_0}(x)\equiv0,\ \ \text{and}\ \  \bar{V}_{\rho_0}(x)\equiv 0,\ \ x\in\Sigma_{\rho_0}.
\end{equation}
Then using the above equalities (\ref{S26}), we immediately obtain
$$0<\bar{u}(x^-_{\rho_0})=u(x)=0, \ \ 0<\bar{v}(x^-_{\rho_0})=v(x)=0, \ \ x\in \partial\mathbb{R}_+^n,$$
which is impossible. Thus $\rho_0=+\infty$ must hold.

Therefore, our goal is to prove (\ref{S26}). Suppose that $(\ref{S26})$ does not hold, then we deduce that
\begin{equation}\label{S27}
\bar{U}_{\rho_0}(x)>0,\ \ \text{and}\ \  \bar{V}_{\rho_0}(x)>0,\ \ x\in\Sigma_{\rho_0}.
\end{equation}
Otherwise, there exists some point $\tilde{x}\in\Sigma_{\rho_0}$ such that $\bar{U}_{\rho_0}(\tilde{x})=0$. We have
$$(-\Delta)^s\bar{U}_{\rho_0}(\tilde{x})=C\int_{\mathbb{R}^n}\frac{-\bar{U}_{\rho_0}(y)}{|\tilde{x}-y|^{n+2s}}dy<0.$$
On the other hand, it is easy to get that
$$(-\Delta)^{s}\bar{U}_\rho(x)=|\tilde{x}^-_\rho|^{-\bar{\tau}_1} |\bar{v}(\tilde{x}^-_\rho)|^{p-1}\bar{v}(\tilde{x}^-_\rho)-|\tilde{x}|^{-\bar{\tau}_1} |\bar{v}(\tilde{x})|^{p-1}\bar{v}(\tilde{x})\geq |\tilde{x}|^{-\bar{\tau}_1} p|\bar{v}|^{p-1}(\tilde{x})\bar{V}_\rho(\tilde{x})\geq 0,$$
where we use the fact $\bar{V}_{\rho_0}\geq 0$. This leads to a contradiction. Hence, (\ref{S27}) holds.

Now we show that  the plane $T_{\rho_0}$ can be moved upward a little bit further and hence obtain a contradiction with the definition of $\rho_0$. Precisely, we will verify that there exists some small $\varepsilon>0$ such that for any $\rho\in (\rho_0,\rho_0+\varepsilon)$,
\begin{equation}\label{S28}
\bar{U}_\rho(x)\geq 0\ \ \text{and}\ \ \bar{V}_\rho(x)\geq 0, \ \ x\in\Sigma_\rho,
\end{equation}
where $\varepsilon$ is determined later.

If (\ref{S28}) is not true, then for any $\varepsilon_k\to 0$ as $k\to+\infty$, there exists $\rho_k\in (\rho_0,\rho_0+\varepsilon_k) $ and $x_k\in\mathbb{R}_+^n\cap\Sigma_{\rho_k}$ such that
\begin{equation}\label{202}
\bar{U}_{\rho_k}(x_k)=\inf_{\Sigma_{\rho_k}}\{\bar{U}_{\rho_k}(x),\bar{V}_{\rho_k}(x)\}<0.
\end{equation}
Similar argument as that of (\ref{201}) gives that
\begin{equation}\label{203}
(-\Delta)^{s_1}\bar{U}_{\rho_k}(x_k)+c(x_k)\bar{U}_{\rho_k}(x_k)\geq 0,\end{equation}
where $c(x)=-|x|^{-\bar{\tau}_1} p\bar{v}^{p-1}(x)$. From (\ref{030}) and the decay conditions of $u$ and $v$, we deduce that
\begin{equation}\label{204}
\mathop{\lim}\limits_{ x \to \infty}|x|^{2s_1}c(x)=0 \text{ and } c(x) \text{ is bounded below in } \Sigma_{\rho_k},
\end{equation}
where we used the assumption $\alpha>-2s_2$. Then from Proposition \ref{L5} we know that there exists $\ell_k>0$ and $R_0>0$ such that
\begin{equation}\label{99}
x_k\in B_{R_0}(0)\cap\Sigma_{\rho_k-\ell_k}.
\end{equation}

Denote $\ell_0 (>0)$ as the constant given in Proposition \ref{L5} corresponding to the half space $\Sigma_{\rho_0+1}$.
Combining the remark about the monotonicity of $\ell$ with respect to $\lambda$ below  Proposition \ref{L5}, (\ref{99}) and the fact that $\varepsilon_k\to 0$, we have that
\begin{equation}\label{206}
x_k\in B_{R_0}(0)\cap \Sigma_{\rho_0-\frac{\ell_0}{2}}.
\end{equation}
If $\rho_0-\frac{\ell_0}{2}\leq 0$, then (\ref{206}) contradicts with the fact that $x_k\in \mathbb{R}^n_+$. If $\rho_0-\frac{\ell_0}{2}>0$, due to (\ref{S27}) and continuity of $\bar{u}$, we know that there exists $\varepsilon'\in(0,\frac{\ell_0}{2})$ such that for any $\varepsilon_k\leq\varepsilon'$ and $\rho\in(\rho_0,\rho_0+\varepsilon_k)$,
$$\bar{U}_{\rho}(x)\geq 0,\ \ \ \  x\in \overline{B_{R_0}(0)\cap \Sigma_{\rho_0-\frac{\ell_0}{2}}}.$$
This contradicts with  (\ref{206}) and (\ref{202}). Hence, we derive that for any $\rho\in (\rho_0,\rho_0+\varepsilon')$ with $\varepsilon'>0$ small enough,
$$\bar{U}_\rho(x)\geq 0,\ \ x\in \Sigma_\rho.$$

Similarly, we may verify that there exists $\varepsilon''>0$ such that for any $\rho\in (\rho_0,\rho_0+\varepsilon'')$ the inequality holds
$$\bar{V}_\rho(x)\geq 0,\ \ x\in \Sigma_\rho.$$
Let $\varepsilon=\min\{\varepsilon',\varepsilon''\}$, then (\ref{S28}) follows immediately and hence (\ref{S26}) holds, which yields that $\rho_0=+\infty$.

The result $\rho_0=+\infty$ indicates that both $\bar{u}(x)$ and  $\bar{v}(x)$ are monotone increasing along the $x_n$-axis. This contradicts with the asymptotic behaviors (\ref{030}). Therefore, $(\bar{u},\bar{v})=(0,0)$, which  yields that $(u,v)=(0,0)$.
We complete the proof of Theorem \ref{T2}.
\endproof

{\bf Acknowledgment.} The  authors are supported by the Natural Science Foundation of Hunan Province, China
(Grant No. 2022JJ30118).

\end{document}